\documentclass[12pt]{article}
\usepackage{amsmath}
\usepackage{amssymb}
\usepackage{amsthm}
\usepackage{tikz}
\usepackage{algorithmic}
\usepackage{algorithm}
\usepackage{url}
\usepackage{multirow}

\newtheorem{theorem}{Theorem}

\newtheorem{conjecture}{Conjecture}
\newtheorem{prop}{Proposition}

\def\sqr#1#2{{\vbox{\hrule height.#2pt
    \hbox{\vrule width.#2pt height#1pt \kern#1pt
        \vrule width.#2pt}\hrule height.#2pt}}}
\def\eqed{\sqr53}
\def\qed{%
    \ifmmode\eqno\eqed
    \else\nobreak\ \hfill\eqed\medbreak\fi}

\newcommand{\rk}{r}
\newcommand{\cl}{{\rm cl}}

\title{Matroids with nine elements}
\author{Dillon Mayhew\\ School of Mathematics, Statistics \& Computer Science\\ Victoria University\\ Wellington, New Zealand\\ 
{\tt Dillon.Mayhew@mcs.vuw.ac.nz} \\ \null\\ \\ 
Gordon F. Royle\\
School of Computer Science \& Software Engineering\\
University of Western Australia\\
Perth, Australia\\
{\tt gordon@csse.uwa.edu.au}
}

\begin{document}
\maketitle

\begin{abstract}
We describe the computation of a catalogue containing all matroids with up to nine elements, and
present some fundamental data arising from this cataogue. Our computation confirms and extends the results obtained in the 1960s by Blackburn, Crapo \& Higgs. The matroids and associated data are stored in an online database, and we give three short examples of the use of this database.

\end{abstract}

\section{Introduction}

In the late 1960s, Blackburn, Crapo \& Higgs published a technical report describing
the results of a computer search for all simple matroids on up to eight elements (although the resulting
paper \cite{MR0419270}  did not appear until 1973). In both the report and the paper they 
said 

\begin{quote}
{\em ``It is unlikely that a complete tabulation of 9-point geometries will be
either feasible or desirable, as there will be many thousands of them. The
recursion $g(9) = g(8)^{3/2}$ predicts 29260."}
\end{quote}

Perhaps this comment dissuaded later researchers in matroid theory, because their catalogue remained unextended for more than 30 years, which surely makes it one of the longest standing 
computational results in combinatorics. However, in this paper we demonstrate that they were in 
fact unduly pessimistic, and describe an orderly algorithm
(see McKay \cite{MR1606516} and Royle \cite{MR1614301}) 
that confirms their computations and extends them by determining the 383172 pairwise 
non-isomorphic matroids on nine elements (see Table~\ref{tab:allmat9}). 

Although this number of matroids is easily manageable on today's 
computers, our  experiments with 10-element matroids suggests that there are at least 
$2.5 \times 10^{12}$ sparse paving matroids of rank 5 on 10 elements. However we refrain from making
any analogous predictions about the desirability or feasibility of constructing a catalogue of 
10-element matroids!

We give some fundamental data about these matroids, and briefly describe how they are
incorporated into an online database that provides access to a far greater range of data; this online
database is accessible at  \url{http://people.csse.uwa.edu.au/gordon/small-matroids.html}.

\section{Matroids, flats and hyperplanes}

We assume that the reader is familiar with the general idea of a matroid as being a 
combinatorial generalization of a multiset of vectors in a vector space, in that it consists 
of a set of elements and combinatorially defined concepts of dependence, 
independence etc. with properties
analogous to the same concepts in vector spaces (see Oxley \cite{MR2065730} for a gentle
introduction to matroids, and Oxley \cite{MR1207587} for complete details).

There are many equivalent ways to make this description precise, but for definiteness we
will take a matroid $M = (E,\rk)$ to be a set $E$ together with a {\rm rank function} 
$$
\rk : 2^E \rightarrow {\mathbb Z}
$$
satisfying the three conditions:

\noindent
(R1) If $A \subseteq E$ then $0 \leq \rk(A) \leq |A|$.\\
(R2) If $A \subseteq B \subseteq E$ then $\rk(A) \leq \rk(B)$.\\
(R3) If $A, B \subseteq E$, then 
$
\rk(A \cap B) + \rk(A \cup B) \leq \rk(A) + \rk(B).
$

The matroid equivalents of various concepts from linear algebra may then be defined using the rank 
function so that, for example, a set of elements $A \subseteq E$ is {\em independent} if 
$|A| = \rk(A)$, and {\em spanning} if $\rk(A) = \rk(E)$. The rank of the matroid itself is 
defined to be $\rk(E)$.

A {\em flat} is a set of elements $F \subseteq E$ such that for all $e \in E \backslash F$,
$$
\rk(F \cup e) > \rk(F).
$$
A flat is the matroid equivalent of a subspace, and as with subspaces, a crucial property of flats is that the intersection of two flats is again a flat. The intersection of all the flats containing a 
subset $A \subseteq E$ is called the {\em closure} of $A$ and denoted by $\cl(A)$; it can be viewed as analogous to the subspace generated by a set of vectors. 

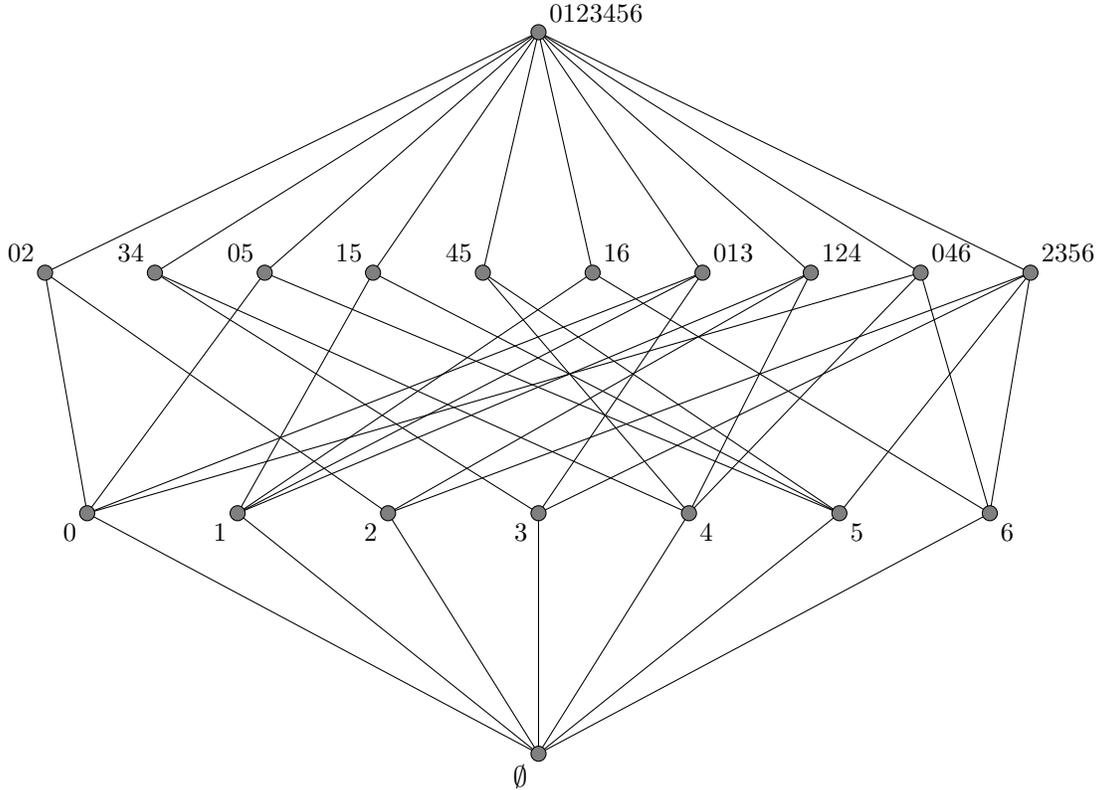
\begin{figure}
\begin{center}
\begin{tikzpicture}[xscale=2, yscale=2]
\draw (4.0,1.6) -- (1.0,3.2);
\draw (4.0,1.6) -- (2.0,3.2);
\draw (4.0,1.6) -- (3.0,3.2);
\draw (4.0,1.6) -- (4.0,3.2);
\draw (4.0,1.6) -- (5.0,3.2);
\draw (4.0,1.6) -- (6.0,3.2);
\draw (4.0,1.6) -- (7.0,3.2);
\draw (1.0,3.2) -- (0.72,4.8);
\draw (1.0,3.2) -- (2.18,4.8);
\draw (1.0,3.2) -- (5.09,4.8);
\draw (1.0,3.2) -- (6.54,4.8);
\draw (2.0,3.2) -- (2.9,4.8);
\draw (2.0,3.2) -- (4.36,4.8);
\draw (2.0,3.2) -- (5.09,4.8);
\draw (2.0,3.2) -- (5.81,4.8);
\draw (3.0,3.2) -- (0.72,4.8);
\draw (3.0,3.2) -- (5.81,4.8);
\draw (3.0,3.2) -- (7.27,4.8);
\draw (4.0,3.2) -- (1.45,4.8);
\draw (4.0,3.2) -- (5.09,4.8);
\draw (4.0,3.2) -- (7.27,4.8);
\draw (5.0,3.2) -- (1.45,4.8);
\draw (5.0,3.2) -- (3.63,4.8);
\draw (5.0,3.2) -- (5.81,4.8);
\draw (5.0,3.2) -- (6.54,4.8);
\draw (6.0,3.2) -- (2.18,4.8);
\draw (6.0,3.2) -- (2.9,4.8);
\draw (6.0,3.2) -- (3.63,4.8);
\draw (6.0,3.2) -- (7.27,4.8);
\draw (7.0,3.2) -- (4.36,4.8);
\draw (7.0,3.2) -- (6.54,4.8);
\draw (7.0,3.2) -- (7.27,4.8);
\draw (0.72,4.8) -- (4.0,6.4);
\draw (1.45,4.8) -- (4.0,6.4);
\draw (2.18,4.8) -- (4.0,6.4);
\draw (2.9,4.8) -- (4.0,6.4);
\draw (3.63,4.8) -- (4.0,6.4);
\draw (4.36,4.8) -- (4.0,6.4);
\draw (5.09,4.8) -- (4.0,6.4);
\draw (5.81,4.8) -- (4.0,6.4);
\draw (6.54,4.8) -- (4.0,6.4);
\draw (7.27,4.8) -- (4.0,6.4);
\filldraw[fill=gray] (4.0,1.6) circle (0.05cm);
\draw(4.0,1.6) node [below left] {\small $\emptyset$};
\filldraw[fill=gray] (1.0,3.2) circle (0.05cm);
\draw(1.0,3.2) node [below left] {\footnotesize $0$};
\filldraw[fill=gray] (2.0,3.2) circle (0.05cm);
\draw(2.0,3.2) node [below left] {\footnotesize $1$};
\filldraw[fill=gray] (3.0,3.2) circle (0.05cm);
\draw(3.0,3.2) node [below left] {\footnotesize $2$};
\filldraw[fill=gray] (4.0,3.2) circle (0.05cm);
\draw(4.0,3.2) node [below left] {\footnotesize $3$};
\filldraw[fill=gray] (5.0,3.2) circle (0.05cm);
\draw(5.0,3.2) node [below right] {\footnotesize $4$};
\filldraw[fill=gray] (6.0,3.2) circle (0.05cm);
\draw(6.0,3.2) node [below right] {\footnotesize $5$};
\filldraw[fill=gray] (7.0,3.2) circle (0.05cm);
\draw(7.0,3.2) node [below right] {\footnotesize $6$};
\filldraw[fill=gray] (0.72,4.8) circle (0.05cm);
\draw(0.72,4.8) node [above left] {\footnotesize $02$};
\filldraw[fill=gray] (1.45,4.8) circle (0.05cm);
\draw(1.45,4.8) node [above left] {\footnotesize $34$};
\filldraw[fill=gray] (2.18,4.8) circle (0.05cm);
\draw(2.18,4.8) node [above left] {\footnotesize $05$};
\filldraw[fill=gray] (2.9,4.8) circle (0.05cm);
\draw(2.9,4.8) node [above left] {\footnotesize $15$};
\filldraw[fill=gray] (3.63,4.8) circle (0.05cm);
\draw(3.63,4.8) node [above left] {\footnotesize $45$};
\filldraw[fill=gray] (4.36,4.8) circle (0.05cm);
\draw(4.36,4.8) node [above right] {\footnotesize $16$};
\filldraw[fill=gray] (5.09,4.8) circle (0.05cm);
\draw(5.09,4.8) node [above right] {\footnotesize $013$};
\filldraw[fill=gray] (5.81,4.8) circle (0.05cm);
\draw(5.81,4.8) node [above right] {\footnotesize $124$};
\filldraw[fill=gray] (6.54,4.8) circle (0.05cm);
\draw(6.54,4.8) node [above right] {\footnotesize $046$};
\filldraw[fill=gray] (7.27,4.8) circle (0.05cm);
\draw(7.27,4.8) node [above right] {\footnotesize $2356$};
\filldraw[fill=gray] (4.0,6.4) circle (0.05cm);
\draw(4.0,6.4) node [above right] {\footnotesize $0123456$};
\end{tikzpicture}
\end{center}
\caption{The lattice of flats of a rank 3 matroid on 7 elements}
\label{fig:lattice}
\end{figure}

The collection of flats of a matroid $M$ under inclusion, together with the two operations 
$$
A \wedge B = A \cap B \qquad A \vee B = \cl(A \cup B)
$$
forms a {\em geometric lattice}  ${\cal L}(M)$. Figure~\ref{fig:lattice} shows this lattice for a rank-3
matroid with 7 elements, where a label such as $013$ is shorthand for the subset $\{0,1,3\}$. 
In this lattice we say that $A$ {\em covers} $B$ if $A \subset B$ and  there are no flats $C$ such 
that $A \subset C \subset B$.
The rank function of $M$ can be 
uniquely recovered from ${\cal L}(M)$ --- the rank of a flat is the length of the maximal
chain from $\emptyset$ to that flat, and the  rank of any other subset $A \subseteq E$ is the rank
of $\cl(A)$. In fact more is true, because the {\em hyperplanes} (i.e. the flats of rank $\rk(E)-1$) are sufficient to determine the remaining flats, and hence the entire matroid.

Two matroids $M_1 = (E_1, \rk_1)$ and $M_2 = (E_2, \rk_2)$ are {\em isomorphic} if there is a 
bijection $\rho: E_1 \rightarrow E_2$ such that $\rk_2(\rho(A)) = \rk_1(A)$ for all $A \subseteq E$.
For most (but not all) applications, it is appropriate to treat isomorphic matroids as equal and when
counting and cataloguing matroids we are usually interested only in pairwise non-isomorphic matroids.

We can determine matroid isomorphism by using 
the {\em hyperplane graph} of the matroid, which is
the bipartite graph whose vertices are the elements and hyperplanes of the matroid and
where a hyperplane-vertex is adjacent to an element-vertex if and only if the hyperplane contains the element. From our previous discussion, it is 
clear that two matroids are isomorphic if and only if their hyperplane graphs are isomorphic as
bipartite graphs (i.e. with the bipartition fixed).
Although the theoretical complexity of graph isomorphism is not known, in practice Brendan
McKay's program {\tt nauty} \cite{MR635936} can easily process graphs with thousands of vertices 
(except for a few pathologically difficult, but poorly understood graphs). As
the hyperplane graphs of the nine-element matroids have an average of only 74 vertices, isomorphism
for matroids of this size is very easy to resolve in practice.

We note that using hyperplanes is a somewhat arbitrary choice and that any other collection 
of subsets that determines the matroid, such as the set of flats or the set of independent sets, could be 
used analogously. 

\section{Extensions and modular cuts}

If $M = (E, \rk)$ is a matroid and $e \in E$, then the restriction of $\rk$ to the subsets of 
$E \backslash e$ is itself a rank function, and so determines a matroid $M\backslash e = 
(E \backslash e, \rk \mid_{E \backslash e})$. We say that $M\backslash e$ is obtained
by {\em deleting} $e$ from $M$ and conversely that $M$ is an {\em single-element extension} of 
$M\backslash e$.

Now suppose that we have a list ${\cal M}_k$ of the matroids on $k$ elements (or more precisely, one representative from each isomorphism class of matroids on $k$ elements). Then we can form 
the list ${\cal M}_{k+1}$ of all matroids on $k+1$ elements by first finding all possible single-element extensions of every matroid in ${\cal M}_k$ and then eliminating unwanted isomorphic copies.

The key to extending a matroid in all possible ways lies in understanding the relationship 
between the flats of a matroid $M$ and the flats of a single-element deletion $N = M\backslash e$. Letting
${\cal F}(M)$ denote the set of flats of a matroid $M$, it is easy to see that
\begin{equation}\label{flats}
{\cal F}(M\backslash e) = \{ F \backslash e \mid F \in {\cal F}(M) \}. 
\end{equation}

Thus suppose that we are given the matroid $N$ and wish to add a new element $e$, thereby finding all matroids $M$ such that $M \backslash e = N$. By (\ref{flats}), every flat of $M$ is of the 
form $F$ or $F \cup \{e\}$ where $F \in {\cal F}(N)$. More precisely, for each flat $F \in {\cal F}(N)$
exactly one of the following three situations must hold  in $M$:

\begin{enumerate}
\item $F \in {\cal F}(M)$ and $F \cup \{e\} \in {\cal F}(M)$
\item $F \in {\cal F}(M)$ but $F \cup \{e\} \notin {\cal F}(M)$
\item $F \notin {\cal F}(M)$ but $F \cup \{e\} \in {\cal F}(M)$
\end{enumerate}

Thus the flats of $M\backslash e$ are partitioned into three parts  in such a way that $M$ can 
be uniquely recovered from $M \backslash e$ and this partition. Thus we can construct every 
possible single-element extension of a matroid $N$ by considering all ``suitable'' partitions of
${\cal F}(N)$ into three parts and forming the different candidates for $M$ accordingly.

This is feasible in practice because Crapo \cite{MR0190045} showed that only certain highly
structured partitions of ${\cal F}(M\backslash e)$ can actually arise, and therefore only a
very limited number of partitions need be considered when extending a matroid.  To describe 
this result we need one more piece of terminology: two flats $F, G \in {\cal F}(M)$ are a {\em
modular pair} if 
$$
\rk(F) + \rk(G) = \rk(F \cup G) + \rk(F \cap G).
$$

What Crapo showed was that if $N = M\backslash e$ and ${\cal F}(N) = {\cal F}_1 \cup {\cal F}_2 \cup {\cal F}_3$ is the partition of the flats of $N$ according to the three possibilities listed above (respectively), then 

\begin{enumerate}
\item ${\cal F}_3$ is an up-set in the lattice ${\cal L}(N)$ i.e. if $F \in {\cal F}_3$ then any flat
containing $F$ is in ${\cal F}_3$.

\item ${\cal F}_3$ is closed under taking intersections of {\em modular pairs} of flats.

\item ${\cal F}_2$ is the set of flats covered in ${\cal L}(N)$ by a member of ${\cal F}_3$.
\end{enumerate}

The set ${\cal F}_3$ is called a {\em modular cut} and ${\cal F}_2$ the {\em collar} of the
modular cut. As the modular cut determines its collar and ${\cal F}_1$ consists of the remaining
flats, it follows that the 
modular cut alone determines the entire partition. Therefore we have the following
result: 

\begin{theorem}
There is a 1-1 correspondence between modular cuts of $N$ and single-element extensions
of $N$. 
\end{theorem}

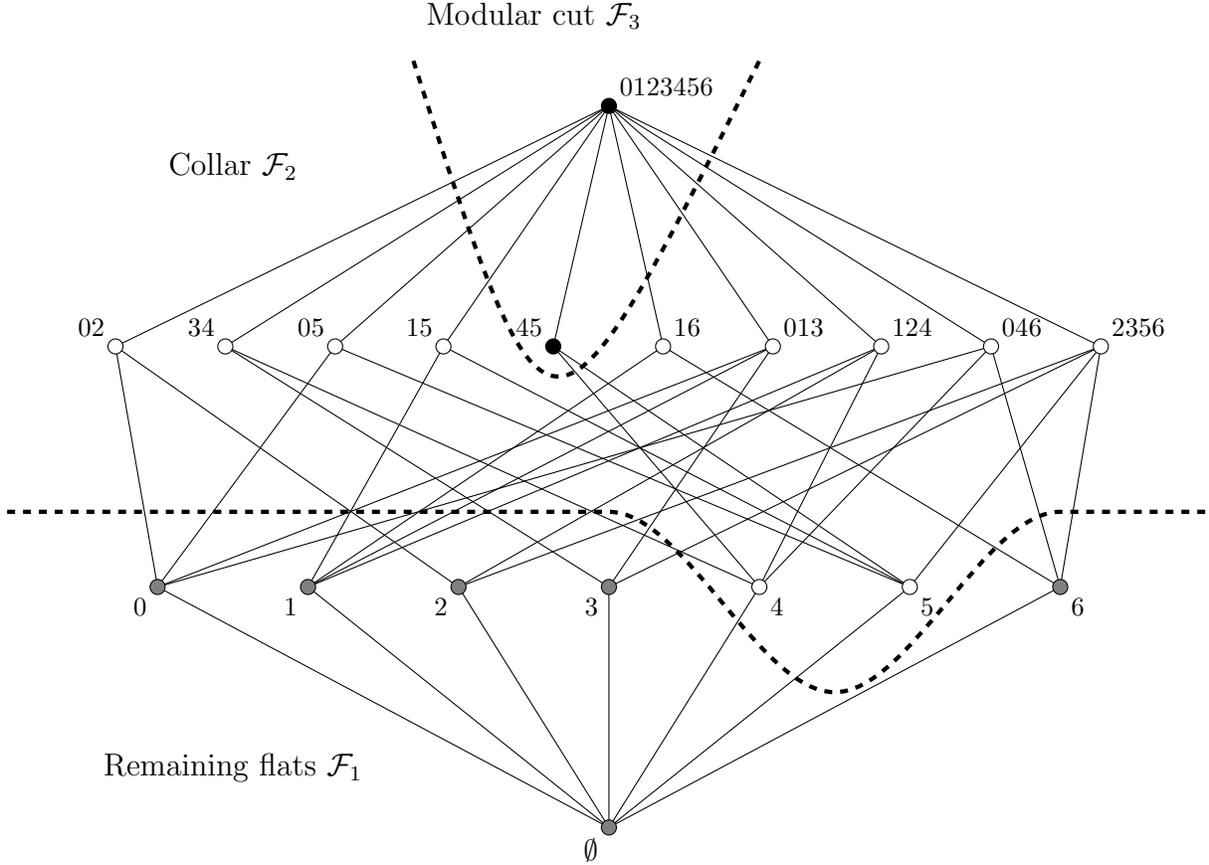
\begin{figure}[t]
\begin{center}
\begin{tikzpicture}[xscale=2, yscale=2]
\draw (4.0,1.6) -- (1.0,3.2);
\draw (4.0,1.6) -- (2.0,3.2);
\draw (4.0,1.6) -- (3.0,3.2);
\draw (4.0,1.6) -- (4.0,3.2);
\draw (4.0,1.6) -- (5.0,3.2);
\draw (4.0,1.6) -- (6.0,3.2);
\draw (4.0,1.6) -- (7.0,3.2);
\draw (1.0,3.2) -- (0.72,4.8);
\draw (1.0,3.2) -- (2.18,4.8);
\draw (1.0,3.2) -- (5.09,4.8);
\draw (1.0,3.2) -- (6.54,4.8);
\draw (2.0,3.2) -- (2.9,4.8);
\draw (2.0,3.2) -- (4.36,4.8);
\draw (2.0,3.2) -- (5.09,4.8);
\draw (2.0,3.2) -- (5.81,4.8);
\draw (3.0,3.2) -- (0.72,4.8);
\draw (3.0,3.2) -- (5.81,4.8);
\draw (3.0,3.2) -- (7.27,4.8);
\draw (4.0,3.2) -- (1.45,4.8);
\draw (4.0,3.2) -- (5.09,4.8);
\draw (4.0,3.2) -- (7.27,4.8);
\draw (5.0,3.2) -- (1.45,4.8);
\draw (5.0,3.2) -- (3.63,4.8);
\draw (5.0,3.2) -- (5.81,4.8);
\draw (5.0,3.2) -- (6.54,4.8);
\draw (6.0,3.2) -- (2.18,4.8);
\draw (6.0,3.2) -- (2.9,4.8);
\draw (6.0,3.2) -- (3.63,4.8);
\draw (6.0,3.2) -- (7.27,4.8);
\draw (7.0,3.2) -- (4.36,4.8);
\draw (7.0,3.2) -- (6.54,4.8);
\draw (7.0,3.2) -- (7.27,4.8);
\draw (0.72,4.8) -- (4.0,6.4);
\draw (1.45,4.8) -- (4.0,6.4);
\draw (2.18,4.8) -- (4.0,6.4);
\draw (2.9,4.8) -- (4.0,6.4);
\draw (3.63,4.8) -- (4.0,6.4);
\draw (4.36,4.8) -- (4.0,6.4);
\draw (5.09,4.8) -- (4.0,6.4);
\draw (5.81,4.8) -- (4.0,6.4);
\draw (6.54,4.8) -- (4.0,6.4);
\draw (7.27,4.8) -- (4.0,6.4);
\filldraw[fill=gray] (4.0,1.6) circle (0.05cm);
\draw(4.0,1.6) node [below left] {\small $\emptyset$};
\filldraw[fill=gray] (1.0,3.2) circle (0.05cm);
\draw(1.0,3.2) node [below left] {\footnotesize $0$};
\filldraw[fill=gray] (2.0,3.2) circle (0.05cm);
\draw(2.0,3.2) node [below left] {\footnotesize $1$};
\filldraw[fill=gray] (3.0,3.2) circle (0.05cm);
\draw(3.0,3.2) node [below left] {\footnotesize $2$};
\filldraw[fill=gray] (4.0,3.2) circle (0.05cm);
\draw(4.0,3.2) node [below left] {\footnotesize $3$};
\filldraw[fill=white] (5.0,3.2) circle (0.05cm);
\draw(5.0,3.2) node [below right] {\footnotesize $4$};
\filldraw[fill=white] (6.0,3.2) circle (0.05cm);
\draw(6.0,3.2) node [below right] {\footnotesize $5$};
\filldraw[fill=gray] (7.0,3.2) circle (0.05cm);
\draw(7.0,3.2) node [below right] {\footnotesize $6$};
\filldraw[fill=white] (0.72,4.8) circle (0.05cm);
\draw(0.72,4.8) node [above left] {\footnotesize $02$};
\filldraw[fill=white] (1.45,4.8) circle (0.05cm);
\draw(1.45,4.8) node [above left] {\footnotesize $34$};
\filldraw[fill=white] (2.18,4.8) circle (0.05cm);
\draw(2.18,4.8) node [above left] {\footnotesize $05$};
\filldraw[fill=white] (2.9,4.8) circle (0.05cm);
\draw(2.9,4.8) node [above left] {\footnotesize $15$};
\filldraw[fill=black] (3.63,4.8) circle (0.05cm);
\draw(3.63,4.8) node [above left] {\footnotesize $45$};
\filldraw[fill=white] (4.36,4.8) circle (0.05cm);
\draw(4.36,4.8) node [above right] {\footnotesize $16$};
\filldraw[fill=white] (5.09,4.8) circle (0.05cm);
\draw(5.09,4.8) node [above right] {\footnotesize $013$};
\filldraw[fill=white] (5.81,4.8) circle (0.05cm);
\draw(5.81,4.8) node [above right] {\footnotesize $124$};
\filldraw[fill=white] (6.54,4.8) circle (0.05cm);
\draw(6.54,4.8) node [above right] {\footnotesize $046$};
\filldraw[fill=white] (7.27,4.8) circle (0.05cm);
\draw(7.27,4.8) node [above right] {\footnotesize $2356$};
\filldraw[fill=black] (4.0,6.4) circle (0.05cm);
\draw(4.0,6.4) node [above right] {\footnotesize $0123456$};
\draw [ultra thick, dashed] (2.7,6.7) .. controls (3.6,3.9) .. (5,6.7);
\draw [ultra thick, dashed] (0,3.7) -- (4,3.7) .. controls +(0:0.5cm) and +(180:0.5cm) .. (5.5,2.5) .. controls +(0:0.5cm) and +(180:0.5cm) .. (7,3.7) -- (8,3.7);
\draw (3.5,7) node {Modular cut ${\cal F}_3$};
\draw (1.5,6) node {Collar ${\cal F}_2$};
\draw (1.5,2) node {Remaining flats ${\cal F}_1$};
\end{tikzpicture}
\end{center}
\caption{A modular cut (black nodes) and its collar (white nodes)}
\label{fig:modcut}
\end{figure}

Figure~\ref{fig:modcut} shows the modular cut $\{45, 0123456\}$ and the corresponding partition for 
the matroid of Figure~\ref{fig:lattice}.

The minimal elements of a modular cut form an anti-chain in ${\cal L}(N)$, and thus an easy way to 
determine the modular cuts of $N$ is simply to compute all the anti-chains of ${\cal L}(N)$, form 
their up-sets and then check that the resulting set of flats is closed under intersection of modular pairs. 

Blackburn, Crapo \& Higgs used a more complicated scheme for computing modular cuts that avoids
creating up-sets that are {\em not} modular cuts, but the overhead of the simpler scheme was 
sufficiently modest that we never had any need to implement the more complicated one. This also
enhances our confidence in the correctness of our results in that existing very-well tested programs (for independent sets in graphs) could be used for computing anti-chains rather than necessarily 
less-tested bespoke programs.

\section{An orderly algorithm}

Our sole remaining task therefore is to consider matroid isomorphism and how to eliminate unwanted
isomorphic copies of the matroids that are constructed, and for this we implemented a 
straightforward (partially) orderly algorithm (Read \cite{MR0491273}, McKay \cite{MR1606516}, Royle \cite{MR1614301}).

In combinatorial construction, an {\em orderly 
algorithm} is one that is structured in such a way that it never outputs more than one
representative of each isomorphism class of the objects being constructed, in this case
matroids.  What this means in practice is that as each matroid is produced by the 
extension procedure, it can be subjected to a test {\em not involving any other matroids} that 
determines whether it should be added to the output or rejected.  Thus there is never 
any need to compare {\em pairs} of matroids, or test a newly-constructed matroid against
a list of previously-constructed ones to check if it is really new.

Our algorithm falls into the category of ``canonical construction path'' orderly algorithms.  Suppose that $M$ is a matroid and that it has hyperplane graph ${\cal H}(M)$. Then {\tt nauty} can be used
to compute the canonical labelling of ${\cal H}(M)$ and thereby identify a distinguished element of
$M$ --- for example, the element that receives the lowest canonical label. This then
identifies a distinguished single-element deletion of $M$, namely the matroid obtained by
deleting the distinguished element. The essence of the canonical construction path orderly
algorithm is that it only accepts matroids that are constructed as an extension of this
distinguished single-element deletion --- whenever an isomorphic copy of $M$ arises as 
an extension of one its {\em other} single-element deletions, it is rejected.  

Although this ensures that the matroids generated by extending one matroid never need
be compared with the matroids generated by extending a different matroid, it is still 
possible that two extensions from the {\em same} matroid may be isomorphic. Indeed this
will necessarily happen if a matroid has two different, but isomorphic, modular cuts. However rather 
than perform isomorph rejection directly on modular cuts (many 
of which may lead to matroids that are subsequently rejected) we instead implemented simple
``compare-and-filter'' isomorph rejection on the set of matroids that were {\em accepted} when 
extending a single matroid.

Putting all this together, we get the procedure described in
Algorithm~\ref{alg1}.

\renewcommand{\algorithmicdo}{{\sc Do}}
\renewcommand{\algorithmicforall}{{\sc For Each}}
\renewcommand{\algorithmicend}{{\sc End}}
\renewcommand{\algorithmicfor}{{\sc For}}

\begin{algorithm}
\caption{Isomorph-free extension of a set $X_k$ of $k$-element matroids}
\label{alg1}
\begin{algorithmic}[1]
\FORALL {matroid $N \in X_k$}
\STATE {Set $N^+ \leftarrow \emptyset$.}
\FORALL {modular cut of $N$}
\STATE {Form the single element extension $M$ determined by the modular cut.}
\STATE {Canonically label ${\cal H}(M)$ and add $M$ to $N^+$ if and only if the newly added element is in the same orbit as the lowest canonically labelled element-vertex of ${\cal H}(M)$.}
\ENDFOR
\STATE {Filter isomorphic matroids from $N^+$ and add the remainder to $X_{k+1}$.}
\ENDFOR
\end{algorithmic}
\end{algorithm}

Notice that each matroid in $X_k$ is processed entirely independently of the remaining matroids
in $X_k$ and therefore the computation can be arbitrarily partitioned between as many 
computers as desired.

\begin{theorem}
If $X_k$ contains one representative from each isomorphism class of $k$-element matroids, then
$X_{k+1}$ contains one representative from each isomorphism class of $(k+1)$-element matroids.
\end{theorem}

\begin{proof}
Let $M$ be an arbitrary $(k+1)$-element matroid, and let $M'$ be its distinguished single-element
deletion. By the hypothesis that $X_k$ contains one
representative from each isomorphism class of $k$-element matroids, a matroid isomorphic to $M'$ will be processed at some stage, and so a matroid isomorphic to $M$ will be constructed and then accepted. The filtering stage ensures that only one isomorph of $M$ will be
accepted during the processing of $M'$ and the orderly aspect of the algorithm ensures that
any isomorph of $M$ is rejected whenever it is constructed as an extension of any matroid other than $M'$. 
Therefore $X_{k+1}$ contains exactly one matroid isomorphic to $M$.
\end{proof}

\section{Results}

We implemented the algorithm described in the previous section, and the resulting numbers
of matroids constructed are summarized in Table~\ref{tab:allmat9} (the totals form 
sequence A055545 in Neil Sloane's OEIS \cite{oeis}).

\begin{table}[h]
\begin{center}
\begin{tabular*}{0.75\textwidth}{@{\extracolsep{\fill}}|l|rrrrrrrrrr|}
\hline
$r\backslash n$& $0$& $1$& $2$& $3$& $4$& $5$& $6$& $7$& $8$ &$9$\\
\hline
$0$&$1$&$1$&$1$&$1$&$1$&$1$&$1$&$1$&$1$&$1$\\
$1$&&$1$&$2$&$3$&$4$&$5$&$6$&$7$&$8$&$9$\\
$2$&&&$1$&$3$&$7$&$13$&$23$&$37$&$58$&$87$\\
$3$&&&&$1$&$4$&$13$&$38$&$108$&$325$&$1275$\\
$4$&&&&&$1$&$5$&$23$&$108$&$940$&$190214$\\
$5$&&&&&&$1$&$6$&$37$&$325$&$190214$\\
$6$&&&&&&&$1$&$7$&$58$&$1275$\\
$7$&&&&&&&&$1$&$8$&$87$\\
$8$&&&&&&&&&$1$&$9$\\
$9$&&&&&&&&&&$1$\\
\hline
Total& $1$& $2$& $4$& $8$& $17$& $38$& $98$& $306$& $1724$&$383172$\\
\hline
\end{tabular*}
\caption{All matroids on up to 9 elements}
\label{tab:allmat9}
\end{center}
\end{table}

These numbers are symmetric with respect to rank because of the theory of matroid {\em duality} --- if
$M = (E, \rk)$ is a matroid of rank $r$, then the function $\rk^*: 2^E \rightarrow {\mathbb Z}$ given
by
$$
\rk^*(A) =  |A| + \rk(E \backslash A) - \rk(E)
$$
is a rank function that determines a matroid $M^* = (E, \rk^*$) of rank $|E|-r$ known as the 
{\em dual} of $M$. We emphasize that our algorithm did {\em not} exploit duality to reduce computation
time by only constructing matroids of rank up to $|E|/2$, but rather we used the fact that our
collection was closed under duality as a ``sanity check'' on the correctness of our implementation.

A {\em loop} in a matroid is an element of rank 0. If a matroid $M$ has a loop $\ell$ then it 
appears in every flat, and the lattice ${\cal L}(M)$ has exactly the same structure as the lattice
${\cal L}(M\backslash \ell)$. Therefore loops play little structural role in $M$ and it is common
to consider them as trivial, and to consider only loopless matroids. 

Two elements $e$ and $f$ are said to be {\em parallel} in a loopless matroid $M$ if $\rk(\{e,f\}) = 1$. In this situation $e$ and $f$ are essentially duplicate elements and each flat either contains both 
of them or neither of them. Replacing $\{e,f\}$ with a single element wherever they occur does not alter the structure of the lattice ${\cal L(M)}$, and so again it is fairly common to consider them as trivial.

A matroid is called {\em simple} if it contains no loops and no parallel elements.  Every matroid 
can be obtained from a unique simple matroid by adding a number of loops and parallel elements,
and so in some sense the simple matroids are the ``building blocks'' from which all matroids
can be constructed.  Conversely, the unique simple matroid obtained from $M$ by removing
all loops and replacing each parallel class (i.e. set of mutually parallel elements) by a single element
is called the {\em simplification} of $M$.

The catalogue of Blackburn, Crapo \& Higgs  only contains the  {\em simple} matroids, and so we give their numbers in Table~\ref{tab:simple9} and note that our computations
are in complete agreement with theirs. In addition, Acketa \cite{MR740223} used Blackburn, Crapo \& Higgs' catalogue to compute the numbers of all matroids (by adding loops and parallel elements
in all possible ways) and our computations are also in complete agreement with his.  More recently, Dukes \cite{MR2119297} has given additional data about the
matroids on up to 8 elements and again our results are in accordance with his.

\begin{table}[h]
\begin{center}
\begin{tabular*}{0.75\textwidth}{@{\extracolsep{\fill}}|l|rrrrrrrrrr|}
\hline
$r\backslash n$& $0$& $1$& $2$& $3$& $4$& $5$& $6$& $7$& $8$ &$9$\\
\hline
$0$&$1$&$$&$$&$$&$$&$$&$$&$$&$$&\\
$1$&&$1$&&&&&&&&\\
$2$&&&$1$&$1$&$1$&$1$&$1$&$1$&$1$&$1$\\
$3$&&&&$1$&$2$&$4$&$9$&$23$&$68$&$383$\\
$4$&&&&&$1$&$3$&$11$&$49$&$617$&$185981$\\
$5$&&&&&&$1$&$4$&$22$&$217$&$188936$\\
$6$&&&&&&&$1$&$5$&$40$&$1092$\\
$7$&&&&&&&&$1$&$6$&$66$\\
$8$&&&&&&&&&$1$&$7$\\
$9$&&&&&&&&&&$1$\\
\hline
Total&$1$&$1$&$1$&$2$&$4$&$9$&$26$&$101$&$950$&$376467$\\
\hline
\end{tabular*}
\caption{Simple matroids on up to 9 elements}
\label{tab:simple9}
\end{center}
\end{table}

We remark that declaring loops and parallel elements --- but not their duals --- to be trivial 
displays a somewhat
{\em graph-theoretical} bias. In a matroid arising from a graph, a loop comes from  
a loop in the graph and parallel elements come from multiple edges, both of which are
routinely excluded in much of graph theory. However, from a matroidal perspective, a matroid and
its dual have equal status, and thus a loop is no more or less trivial than its dual, which is a {\em coloop}. 
Similarly, parallel elements are no more or less trivial than the dual structure, which are elements
in {\em series}. In graphs, coloops correspond to cut-edges and elements in series correspond to paths
with internal vertices of degree two. Graph theorists are understandably reluctant to declare these 
structures trivial because it would mean doing away with both trees and cycles!  Matroidally 
however, the natural building blocks are those matroids that are both {\em simple} and {\em cosimple} (i.e. the dual matroid is also simple) and so we give their numbers in Table~\ref{tab:simplecosimple9}.

\begin{table}[h]
\begin{center}
\begin{tabular*}{0.75\textwidth}{@{\extracolsep{\fill}}|l|rrrrrrrrrr|}
\hline
$r\backslash n$& $0$& $1$& $2$& $3$& $4$& $5$& $6$& $7$& $8$ &$9$\\
\hline
$0$&$1$&$$&$$&$$&$$&$$&$$&$$&$$&\\
$1$&&&&&&&&&&\\
$2$&&&&&$1$&$1$&$1$&$1$&$1$&$1$\\
$3$&&&&&&$1$&$6$&$20$&$65$&$380$\\
$4$&&&&&&&$1$&$20$&$525$&$185620$\\
$5$&&&&&&&&$1$&$65$&$185620$\\
$6$&&&&&&&&&$1$&$380$\\
$7$&&&&&&&&&&$1$\\
\hline
Total&$1$&$0$&$0$&$0$&$1$&$2$&$8$&$42$&$657$&$372002$\\
\hline
\end{tabular*}
\caption{Simple and cosimple matroids on up to 9 elements}
\label{tab:simplecosimple9}
\end{center}
\end{table}


\section{Paving Matroids}

A {\em circuit} in a matroid is a minimal dependent set. It is possible for a matroid to have
no circuits (in which case it consists entirely of coloops) but otherwise a matroid of rank $r$
must have a circuit of size at most $r+1$. If the minimum circuit size is equal to $r+1$, then the
matroid is a {\em uniform matroid} $U_{r,n}$ which has the property that the rank of a set $A$
is equal to $\min(r, |A|)$. If the minimum circuit size is at least $r$, then the matroid is called a {\em paving matroid}.  

We need some more terminology before we can understand why paving matroids
form an important class of matroids. 


A {\em $d$-partition} of a set $E$ is a set ${\cal S}$ of subsets of $E$ all of size at least $d$, such that every $d$-subset of $E$ lies in a unique member of ${\cal S}$. Therefore a $1$-partition of a set is simply a normal partition, while a $2$-partition of a set is known as
a {\em pairwise balanced design} with index 1. Obviously the set ${\cal S} = \{E\}$ is a 
$d$-partition for any $d$, and we call this the {\em trivial} $d$-partition. 

The connection between paving matroids and $d$-partitions is given by the following 
result:

\begin{theorem}
If $M = (E,\rk)$ is a paving matroid of rank $d+1 \geq 2$ then its hyperplanes form a non-trivial $d$-partition of $E$. Conversely, the elements of any non-trivial $d$-partition of $E$ form the set of
hyperplanes of a paving matroid of rank $d+1$. 
\end{theorem}

The {\em discrete} $d$-partition of a set $E$ consists of all the $d$-subsets of $E$ and the
corresponding paving matroid is the uniform matroid $U_{d+1,|E|}$.

Based on the rather limited evidence in the catalogue of matroids on up to 8 elements, 
Welsh \cite{MR0427112} asked whether {\em most} matroids are paving matroids. Examining the catalogue of 9-element 
matroids and tabulating the results in Table~\ref{tab:simplepaving9} we see that $71.71\%$ of the simple matroids on 9 elements are paving matroids, compared to $49.50\%$ of the 8-element simple matroids, thus providing some additional evidence that paving matroids do indeed predominate.

\begin{table}[h]
\begin{center}
\begin{tabular*}{0.75\textwidth}{@{\extracolsep{\fill}}|l|rrrrrrrrrr|}
\hline
$r\backslash n$& $0$& $1$& $2$& $3$& $4$& $5$& $6$& $7$& $8$ &$9$\\
\hline
$0$&$1$&$$&$$&$$&$$&$$&$$&$$&$$&\\
$1$&&$1$&&&&&&&&\\
$2$&&&$1$&$1$&$1$&$1$&$1$&$1$&$1$&$1$\\
$3$&&&&$1$&$2$&$4$&$9$&$23$&$68$&$383$\\
$4$&&&&&$1$&$2$&$5$&$18$&$322$&$147163$\\
$5$&&&&&&$1$&$2$&$5$&$39$&$119050$\\
$6$&&&&&&&$1$&$2$&$6$&$178$\\
$7$&&&&&&&&$1$&$2$&$6$\\
$8$&&&&&&&&&$1$&$2$\\
$9$&&&&&&&&&&$1$\\
\hline
Total&$1$&$1$&$1$&$2$&$4$&$8$&$18$&$50$&$439$&$266784$\\
\hline
\end{tabular*}
\caption{Simple paving matroids on up to 9 elements}
\label{tab:simplepaving9}
\end{center}
\end{table}

A $d$-partition is called {\em sparse} if it contains no subsets of size greater than 
$d+1$, and similarly we call a paving matroid of rank $d+1$ sparse if its hyperplanes all have size
$d$ or $d+1$.  A sparse paving matroid is determined completely by its hyperplanes of 
size $d+1$ --- the $d$-partition must consist of these hyperplanes together with every $d$-set
not yet contained in one of these. These hyperplanes are necessarily circuits, and so they
form the set of {\em circuit-hyperplanes} of the matroid.  

Sparse paving matroids have the attractive property that their duals are also sparse paving
matroids --- in fact the circuit-hyperplanes of $M^*$ are the complements of the 
circuit-hyperplanes of $M$. Moreover if $M$ and its dual are both paving matroids, then 
they are necessarily sparse and so the sparse paving matroids forms the largest possible dual-closed
family of paving matroids.  

Computationally, sparse paving matroids are attractive
because they can be viewed simply as independent sets in a certain graph. 
The {\em Johnson graph} $J(n,d+1)$ is the graph whose vertices are all the $(d+1)$-subsets
of an $n$-set, and where two vertices are adjacent if and only if the intersection of the corresponding
subsets has size $d$. Therefore an {\em independent set} of vertices in $J(n,d+1)$ 
is precisely the set of circuit-hyperplanes of a sparse paving matroid of rank $d+1$, 
and conversely. Moreover the automorphism group of $J(n,d+1)$ is equal to the symmetric group 
$S_n$ except in the special case where $n = 2(d+1)$ in which case the graph has an additional automorphism of order 2 induced by complementation on $(d+1$)-sets.

We note in passing that the {\em size} of the maximum independent set in the Johnson graphs $J(n,d+1)$ has been intensively
studied because such an independent set is directly equivalent to a {\em constant weight
code} of length $n$, weight $d+1$ and minimum distance $4$.

More generally, we can form the analogous graph on the $d+1$, $d+2$, $\ldots$, $n-1$ sets of an 
$n$-set where again two vertices are adjacent if the corresponding sets meet in set of size $d$. Then
an independent set in this graph corresponds to a not-necessarily-sparse paving matroid

\section{Representability}

A matroid $M = (E,r)$ of rank $k$ is {\em representable} over a field ${\mathbb F}$ if there is a mapping
$$
\rho: E \rightarrow {\mathbb F}^k
$$
such that for any set $A \subseteq E$,
$$
r(A) = \dim {\rm span}(\rho(A)).
$$

To prove that a matroid {\em is} representable, it suffices to provide a suitable representation $\rho$, but it is considerably harder to prove that a matroid is {\em not} representable. However, Ingleton showed that if $M = (E,r)$ is representable and $A$, $B$, $C$,  $D \subseteq E$, then 
\begin{align*}
r(A) + r(B) + r(A \cup B \cup C) & +  r(A \cup B \cup D) + r(C \cup D) \\
& \leq  r(A \cup B) + r(A \cup C) + r(A \cup D) + r(B \cup C) + r(B \cup D).
\end{align*}
It is therefore sometimes possible to show that a matroid is not representable by displaying four
subsets $A$, $B$, $C$ and $D$ for which this inequality is violated. For want of a convenient term, we will call such matroids {\em Ingleton non-representable}.

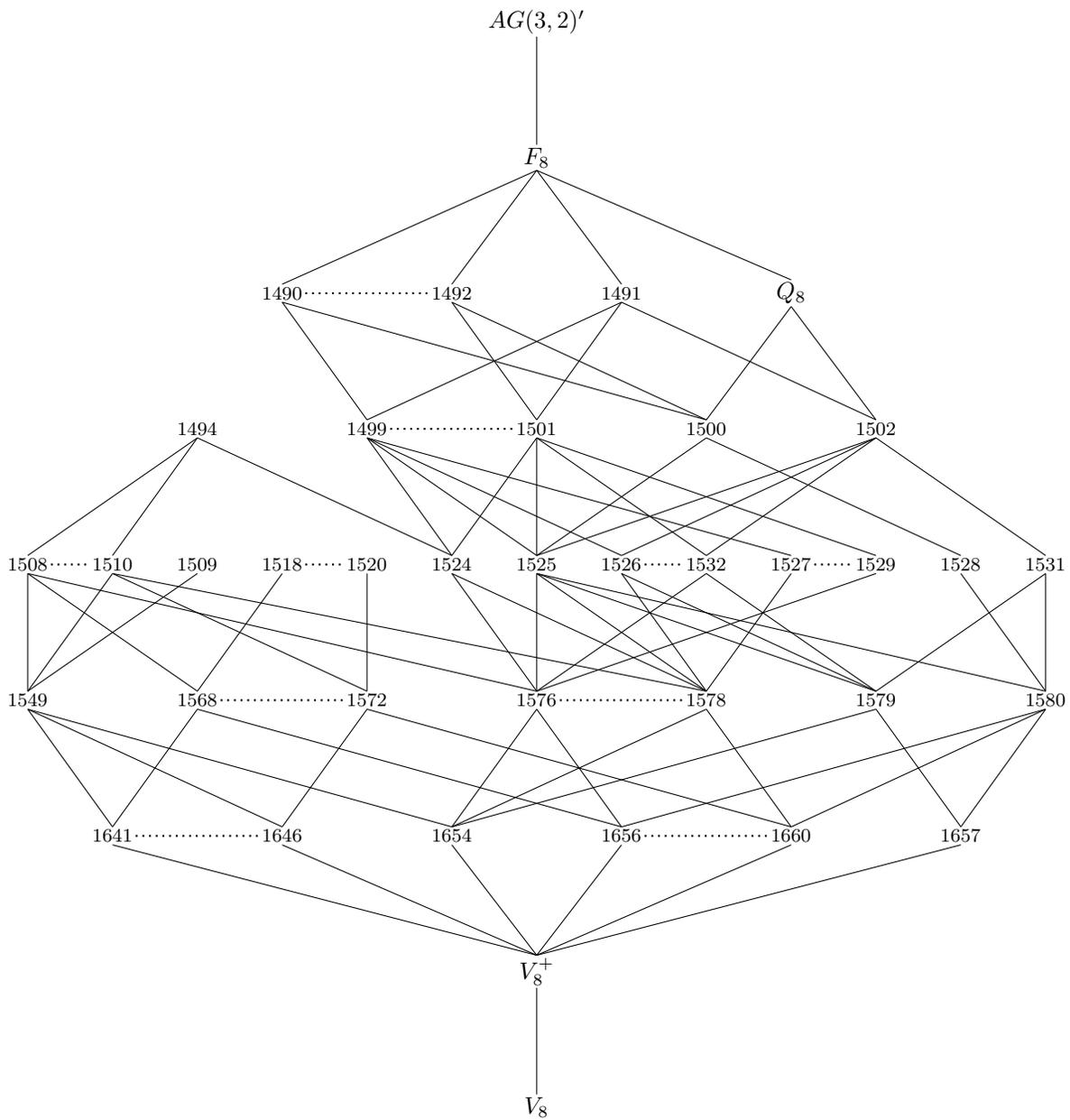
\begin{figure}
\begin{center}
\begin{tikzpicture}[xscale=1.25,yscale =2]

\draw (6,0) node [inner sep = 1pt] (v1771) {{\footnotesize $V_8$}};
\draw (6,1) node [inner sep = 1pt] (v1724) {{\footnotesize $V_8^+$}};

\draw(6,8) node [inner sep = 1pt] (v1485) {{\footnotesize $AG(3,2)'$}};

\draw(6,7) node [inner sep = 1pt] (v1486) {{\footnotesize $F_8$}};

\draw(3,6) node [inner sep = 1pt] (v1490) {{\scriptsize $1490$}};
\draw(5,6) node [inner sep = 1pt] (v1492) {{\scriptsize $1492$}};
\draw(7,6) node [inner sep = 1pt] (v1491) {{\scriptsize $1491$}};
\draw(9,6) node [inner sep = 1pt] (v1493) {{\footnotesize $Q_8$}};

\draw(2,5) node [inner sep = 1pt] (v1494) {{\scriptsize $1494$}};
\draw(4,5) node [inner sep = 1pt] (v1499) {{\scriptsize $1499$}};
\draw(6,5) node [inner sep = 1pt] (v1501) {{\scriptsize $1501$}};
\draw(8,5) node [inner sep = 1pt] (v1500) {{\scriptsize $1500$}};
\draw(10,5) node [inner sep = 1pt] (v1502) {{\scriptsize $1502$}};

\draw(0,4) node [inner sep = 1pt] (v1508) {{\scriptsize $1508$}};
\draw(1,4) node [inner sep = 1pt] (v1510) {{\scriptsize $1510$}};
\draw(2,4) node [inner sep = 1pt] (v1509) {{\scriptsize $1509$}};
\draw(3,4) node [inner sep = 1pt] (v1518) {{\scriptsize $1518$}};
\draw(4,4) node [inner sep = 1pt] (v1520) {{\scriptsize $1520$}};
\draw(5,4) node [inner sep = 1pt] (v1524) {{\scriptsize $1524$}};
\draw(6,4) node [inner sep = 1pt] (v1525) {{\scriptsize $1525$}};
\draw(7,4) node [inner sep = 1pt] (v1526) {{\scriptsize $1526$}};
\draw(8,4) node [inner sep = 1pt] (v1532) {{\scriptsize $1532$}};
\draw(9,4) node [inner sep = 1pt] (v1527) {{\scriptsize $1527$}};
\draw(10,4) node [inner sep = 1pt] (v1529) {{\scriptsize $1529$}};
\draw(11,4) node [inner sep = 1pt] (v1528) {{\scriptsize $1528$}};
\draw(12,4) node [inner sep = 1pt] (v1531) {{\scriptsize $1531$}};

\draw(0,3) node [inner sep = 1pt] (v1549) {{\scriptsize $1549$}};
\draw(2,3) node [inner sep = 1pt] (v1568) {{\scriptsize $1568$}};
\draw(4,3) node [inner sep = 1pt] (v1572) {{\scriptsize $1572$}};
\draw(6,3) node [inner sep = 1pt] (v1576) {{\scriptsize $1576$}};
\draw(8,3) node [inner sep = 1pt] (v1578) {{\scriptsize $1578$}};
\draw(10,3) node [inner sep = 1pt] (v1579) {{\scriptsize $1579$}};
\draw(12,3) node [inner sep = 1pt] (v1580) {{\scriptsize $1580$}};

\draw (1,2) node [inner sep=1pt] (v1641) {{\scriptsize $1641$}};
\draw (3,2) node [inner sep=1pt] (v1646) {{\scriptsize $1646$}};
\draw (5,2) node [inner sep=1pt] (v1654) {{\scriptsize $1654$}};
\draw (7,2) node [inner sep=1pt] (v1656) {{\scriptsize $1656$}};
\draw (9,2) node [inner sep=1pt] (v1660) {{\scriptsize $1660$}};
\draw (11,2) node [inner sep=1pt] (v1657) {{\scriptsize $1657$}};

\draw (v1485.south)--(v1486.north);
\draw (v1486.south)--(v1491.north);
\draw (v1486.south)--(v1490.north);
\draw (v1486.south)--(v1493.north);
\draw (v1486.south)--(v1492.north);
\draw (v1490.south)--(v1499.north);
\draw (v1490.south)--(v1500.north);
\draw (v1491.south)--(v1502.north);
\draw (v1491.south)--(v1501.north);
\draw (v1491.south)--(v1499.north);
\draw (v1492.south)--(v1501.north);
\draw (v1492.south)--(v1500.north);
\draw (v1493.south)--(v1502.north);
\draw (v1493.south)--(v1500.north);
\draw (v1494.south)--(v1510.north);
\draw (v1494.south)--(v1524.north);
\draw (v1494.south)--(v1508.north);
\draw (v1499.south)--(v1527.north);
\draw (v1499.south)--(v1525.north);
\draw (v1499.south)--(v1526.north);
\draw (v1499.south)--(v1524.north);
\draw (v1500.south)--(v1525.north);
\draw (v1500.south)--(v1528.north);
\draw (v1501.south)--(v1532.north);
\draw (v1501.south)--(v1529.north);
\draw (v1501.south)--(v1525.north);
\draw (v1501.south)--(v1524.north);
\draw (v1502.south)--(v1531.north);
\draw (v1502.south)--(v1525.north);
\draw (v1502.south)--(v1532.north);
\draw (v1502.south)--(v1526.north);
\draw (v1508.south)--(v1549.north);
\draw (v1508.south)--(v1576.north);
\draw (v1508.south)--(v1568.north);
\draw (v1509.south)--(v1549.north);
\draw (v1510.south)--(v1572.north);
\draw (v1510.south)--(v1578.north);
\draw (v1510.south)--(v1549.north);
\draw (v1518.south)--(v1568.north);
\draw (v1520.south)--(v1572.north);
\draw (v1524.south)--(v1578.north);
\draw (v1524.south)--(v1576.north);
\draw (v1525.south)--(v1578.north);
\draw (v1525.south)--(v1580.north);
\draw (v1525.south)--(v1579.north);
\draw (v1525.south)--(v1576.north);
\draw (v1526.south)--(v1578.north);
\draw (v1526.south)--(v1579.north);
\draw (v1527.south)--(v1578.north);
\draw (v1528.south)--(v1580.north);
\draw (v1529.south)--(v1576.north);
\draw (v1531.south)--(v1579.north);
\draw (v1531.south)--(v1580.north);
\draw (v1532.south)--(v1576.north);
\draw (v1532.south)--(v1579.north);
\draw (v1549.south)--(v1646.north);
\draw (v1549.south)--(v1654.north);
\draw (v1549.south)--(v1641.north);
\draw (v1568.south)--(v1641.north);
\draw (v1568.south)--(v1656.north);
\draw (v1572.south)--(v1646.north);
\draw (v1572.south)--(v1660.north);
\draw (v1576.south)--(v1654.north);
\draw (v1576.south)--(v1656.north);
\draw (v1578.south)--(v1660.north);
\draw (v1578.south)--(v1654.north);
\draw (v1579.south)--(v1654.north);
\draw (v1579.south)--(v1657.north);
\draw (v1580.south)--(v1660.north);
\draw (v1580.south)--(v1657.north);
\draw (v1580.south)--(v1656.north);
\draw (v1641.south)--(v1724.north);
\draw (v1646.south)--(v1724.north);
\draw (v1654.south)--(v1724.north);
\draw (v1656.south)--(v1724.north);
\draw (v1657.south)--(v1724.north);
\draw (v1660.south)--(v1724.north);
\draw (v1724.south)--(v1771.north);

\draw [thick, dotted] (v1490)--(v1492);
\draw [thick, dotted] (v1499)--(v1501);
\draw [thick, dotted] (v1508)--(v1510);
\draw [thick, dotted] (v1518)--(v1520);
\draw [thick, dotted] (v1526)--(v1532);
\draw [thick, dotted] (v1527)--(v1529);
\draw [thick, dotted] (v1568)--(v1572);
\draw [thick, dotted] (v1576)--(v1578);
\draw [thick, dotted] (v1641)--(v1646);
\draw [thick, dotted] (v1656)--(v1660);

\end{tikzpicture}
\caption{The 39 Ingleton-non-representable matroids on 8 elements}
\label{fig:nonreps}
\end{center}
\end{figure}

Figure~\ref{fig:nonreps} gives a schematic diagram of the Ingleton-non-representable matroids on 8 elements; all of the matroids are sparse paving matroids and the diagram shows how they are related to each other under {\em relaxation of a circuit-hyperplane}, so that for example, the matroid $F_8$ is
obtained from $AG(3,2)'$ by relaxing a circuit-hyperplane. Any named matroids are listed according to their names in Oxley \cite{MR1207587} while the remainder are given just by their number in the database. Dual pairs of matroids are connected by dotted lines.

In addition to the 39 Ingleton-non-representable matroids, there are five other rank-4 matroids on 8
elements that are non-representable. Four of these are related to the sparse paving matroid $P_8$. 
The matroid $P_8$ is a ternary matroid which has the representation

$$
\bordermatrix{&0&1&2&3&\phantom{-}4&5&6&\phantom{-}7\cr
& 1 & 0 & 0 & 0 & \phantom{-}0 & 1 & 1 & -1 \cr
& 0 & 1 & 0 & 0 & \phantom{-}1 & 0 & 1 & \phantom{-}1 \cr
& 0 & 0 & 1 & 0 & \phantom{-}1 & 1 & 0 & \phantom{-}1 \cr
& 0 & 0 & 0 & 1 & -1 & 1 & 1 &\phantom{-} 0\cr
}
$$

The circuit-hyperplanes of this matroid are $\{0,1,2,3\}$, $\{0,1,3,6\}$, $ \{0,2,3,5\}$, $\{1,2,3,4\}$, $\{0,3,4,7\}$,
$\{1,2,5,6\}$, $\{0,4,5,6\}$, $\{1,4,5,7\}$, $\{2,4,6,7\}$ and $\{3,5,6,7\}$. We define four associated sparse paving matroids as follows: $P_1$ is obtained from $P_8$ by relaxing the circuit-hyperplane $\{3,5,6,7\}$, $P_2'$
is obtained from $P_1$ by relaxing $\{0,3,4,7\}$, $P''_2$ is obtained from $P_1$ by relaxing $\{1,2,5,6\}$ and
$P_3$ is obtained from $P_1$ by relaxing both $\{0,3,4,7\}$ and $\{1,2,5,6\}$. 

\begin{prop}
The four matroids $P_1$, $P'_2$, $P''_2$ and $P_3$ are all non-representable, but not Ingleton non-representable.
\end{prop}

\begin{proof}
Let $M \in \{P_1, P'_2, P''_2, P_3\}$ and consider the basis $B = \{0,1,2,3\}$. Then following Section 6.4 of Oxley \cite{MR1207587} a representation for $M$ may be assumed to have the following form where $a,b,c,d,e \not=0$ are
unknown elements of some field.

\[
A = 
\bordermatrix{&0&1&2&3&4&5&6&7\cr
& 1 & 0 & 0 & 0 & 0 & 1 & 1 & 1 \cr
& 0 & 1 & 0 & 0 & 1 & 0 & 1 & a \cr
& 0 & 0 & 1 & 0 & 1 & b & 0 & c \cr
& 0 & 0 & 0 & 1 & 1& d & e & 0\cr
}
\]

Now the sets $\{0,4,5,6\}$, 
$\{1,4,5,7\}$ and $\{2,4,6,7\}$  are circuits in $M$ and so the submatrices of $A$ defined on those particular sets of columns each have determinant 0. This gives the following three conditions respectively: $b(e-1) + d = 0$, 
$b - c - d = 0$ and $a+e -1 = 0$ which implies that
$$
a = (1-e), \quad c = be\quad {\rm and} \quad d = b(1-e)
$$

However consider the submatrix of $A$ with columns $\{3,5,6,7\}$. The determinant of this is
\[
\left|\begin{array}{cccc}0 & 1 & 1 & 1 \\0 & 0 & 1 & 1-e \\0 & b & 0 & be \\1 & b(1-e) & e & 0\end{array}\right| = 0
\]
contradicting the fact that $\{3,5,6,7\}$ is independent in $M$.

Checking Ingleton non-representability of a matroid is a task best left to a computer. \end{proof}

The fifth non-representable matroid of size 8 for which Ingleton's condition gives no information is obtained from the matroid $L_8$ by relaxing a circuit-hyperplane, where $L_8$ is the sparse paving matroid whose circuit-hyperplanes are the 6 faces and the 2 colour-classes of a cube (see Oxley \cite{MR1207587} p510). It can be shown to be non-representable using an analogous argument.


How effective is Ingleton's criterion for detecting non-representability among 9-element matroids? Perhaps surprisingly, it gives {\em no additional information} at all --- that is,  a 9-element matroid is Ingleton non-representable if and only if it contains an Ingleton-non-representable matroid on 8 elements as a minor.

\section{A matroid database}

One of the major uses of any sort of combinatorial catalogue is to compile data regarding the
various combinatorial properties of the objects in the catalogue, and then to use this to
answer questions or explore conjectures concerning the existence, or number of objects
with various combinations of properties. 

A common limitation of combinatorial catalogues is that their use is often restricted to their
immediate creator and/or those researchers willing and able to download the raw data files and write 
their own programs, often resulting in significant duplication of effort. 
We have attempted to ameliorate this problem by incorporating the data into a 
relational database (using MySQL, although this is not important) and providing an online 
interface that permits ``end users" to search, browse and investigate the data. 

Currently we have computed a fairly substantial subset of what might be termed the ``fundamental properties" of
matroids. This includes various counts associated with each matroid such as numbers of loops, coloops, circuits, cocircuits, 
independent sets, bases, hyperplanes, flats and circuit-hyperplanes. It includes numerical properties such as the
size of the automorphism group, the number of orbits of the automorphism group, the connectivity, the minimum circuit size.
Various structural properties such as whether the matroid is binary, ternary, regular, paving, base-orderable, transversal and so on have
also been included.  More importantly however, we have incorporated information about the {\em relationships} between
the matroids --- relationships such as duality, deletion and contraction of elements, relaxation of circuit-hyperplanes, truncations and simplifications. Finally we have included auxiliary information such as information about rank polynomials and representations over small finite fields.

Rather than present a large number of tables of data in this paper, we give three simple examples of the use of the database, and invite readers to explore their own particular interests by using the 
database at \url{http://people.csse.uwa.edu.au/gordon/small-matroids.html}.

\subsection{Excluded minors for $GF(5)$}

One of the most fundamental results in matroid theory is Tutte's characterization of matroids
representable over $GF(2)$ in terms of excluded minors: a matroid is binary if and only if it does
not contain $U_{2,4}$ as a minor. 

Similar characterizations are known for matroids representable over $GF(3)$ where there
are 4 excluded minors and $GF(4)$ where there are 7 excluded minors.  The
analogous characterization for $GF(5)$ is not known or even conjectured, with prevailing
opinion suggesting that such a characterization is likely to be extremely complex and unwieldy. 
In fact, Whittle \cite{MR2179649}
suggests that ``It is not clear that the problem for finding the specific excluded minors
for $GF(5)$ is that well motivated'' and that the real question in representability is
to resolve Rota's conjecture that the list of excluded minors for representability over any
field is finite.

It is straightforward to determine the matroids on up to 9 elements that {\em are} representable over $GF(5)$ by finding all sets of at most 9 points in the projective space $PG(3,5)$ that are pairwise inequivalent under the action of the  group $PGL(4,5)$, identifying the corresponding matroids (which have rank at most 4) and then finding their duals.

Once this is done, we can identify the matroids that are {\em not} $GF(5)$-representable but for
which every single-element deletion and single-element contraction {\em is} $GF(5)$-representable
and thus determine the excluded minors on at most 9 elements. Table~\ref{tab:exgf5} shows
the numbers of matroids that were found, confirming the belief that an excluded minor
characterization of $GF(5)$-representable matroids in the traditional style would indeed be very cumbersome.

\begin{table}[ht]
\begin{center}
\begin{tabular}{|ccrl|}
\hline
Size & Rank & No.&Comment\\
\hline
7&2&1&Uniform $U_{2,7}$\\
7&3&5&\\
7&4&5&\\
7&5&1&Uniform $U_{5,7}$\\
\hline
8&3&2&\\
8&4&92&\\
8&5&2&\\
\hline
9&3&9&\\
9&4&219&\\
9&5&219&\\
9&6&9&\\
\hline
\end{tabular}
\end{center}
\caption{Excluded minors for $GF(5)$ on up to nine elements}
\label{tab:exgf5}
\end{table}

\subsection{Numbers of bases}

In a matroid of rank $r$ on $n$ elements, the number $b$ of bases must necessarily satisfy
$
1 \leq b \leq \binom{n}{r}.
$
In 1969, Welsh \cite{MR0278975} conjectured that for {\em every triple} $(n,r,b)$ such that $0 \leq r \leq n$ and $1 \leq b \leq \binom{n}{r}$, there is
a matroid of rank $r$ on $n$ elements with exactly $b$ bases --- in other words, everything that can happen, does.

We can check this all matroids on up to 9 elements with a single SQL statement (though note that \verb+binomial+ is not a built-in function, but must be programmed):

{\small
\begin{verbatim}
SELECT tmp.size, tmp.rank, COUNT(*) FROM 
  (SELECT DISTINCT size, rank, numBases FROM matroids9) as tmp
  GROUP BY tmp.size, tmp.rank 
  HAVING COUNT(*) <> binomial(tmp.size, tmp.rank);
\end{verbatim}
}

The inner \verb+SELECT+ statement first creates a list of all the distinct triples $(n,r,b)$ represented in the database and
gives it the alias \verb+tmp+. The outer \verb+SELECT .... GROUP BY+ statement {\em counts} the triples in \verb+tmp+ for 
each fixed pair $(n,r)$, while the \verb+HAVING+ statement extracts the pairs where this count is not equal to $\binom{n}{r}$, thus representing one or
more ``missing'' triples. 

{\small
\begin{verbatim}
+------+------+----------+
| size | rank | COUNT(*) |
+------+------+----------+
|    6 |    3 |       19 | 
+------+------+----------+
\end{verbatim}
}

The output  shows that there are only 19 triples of the form $(6,3,b)$, rather than the expected 20. In fact, the missing triple is $(6,3,11)$ --- there are no rank-3 matroids on 6 elements with exactly 11 bases --- a fact that was previously observed by Anna de Mier (personal communication). The absence of any  other missing triples with $n \leq 9$ and the exponential explosion in numbers of matroids as $n$ reaches 10 leads us to strongly believe the following conjecture:

\begin{conjecture}
For every triple $(n,r,b)$ such that $0 \leq r \leq n$ and $1 \leq b \leq \binom{n}{r}$ there is
a matroid of rank $r$ on $n$ elements with exactly $b$ bases {\em except when}
$$
(n,r,b) = (6,3,11).
$$
\end{conjecture}

\subsection{Transversal matroids}

Given two bases $A$ and $B$ of a matroid, the subsets $X \subseteq A$ and
$Y \subseteq B$ are called {\em exchangeable} if both 
$A \backslash X \cup Y$ and $B \backslash Y \cup X$
are bases.

A matroid is called {\em base-orderable} if for any two bases $A$ and $B$ there is 
a bijection $\varphi: A \rightarrow B$ such that $a$ and $\varphi(a)$ are exchangeable, while it
is {\em strongly base-orderable} if the bijection can be selected so that $X$ and $\varphi(X)$
are exchangeable for every {\em subset} $X \subseteq A$. A matroid on a set $E$ is a {\em transversal matroid} if there is a bipartite graph $G$ with bipartition $E \cup F$ such that the independent sets of $M$ are precisely the subsets of $E$ that are the endpoints of a matching (i.e. an independent set of edges) of $G$. 

It is obvious that strongly base-orderable matroids are base-orderable, but less obvious that 
transversal matroids and their duals (which need not be transversal matroids) are strongly-base orderable. These classes of matroids are important but not fully understood, and 
therefore the numbers of  matroids in each of these classes is of some interest. Determining whether a matroid is base-orderable or strongly base-orderable is straightforward, and we implemented the 
algorithm given by Brualdi \& Dinolt \cite{MR0304210} for testing transversality. 

Table~\ref{tab:transversal} gives these
numbers where each cell of the table contains four numbers which, reading from top to
bottom are the total number of matroids and the number of base-orderable, strongly base-orderable and 
transversal matroids respectively. For the omitted ranks (ranks 0, 1, 7, 8 and 9) all the matroids
in the catalogue are
transversal.

\begin{table}[t]
\begin{center}
\begin{tabular}{|c|rrrrrrrr|}
\hline
Rank $\backslash$Size & $2$ &$3$&$4$&$5$&$6$&$7$&$8$&$9$\\
\hline
\multirow{4}{*}{$2$}&\small $1$&\small $3$&\small $7$&\small $13$&\small $23$&\small $37$&\small $58$&\small $87$\\
&\small $1$&\small $3$&\small $7$&\small $13$&\small $23$&\small $37$&\small $58$&\small $87$\\
&\small $1$&\small $3$&\small $7$&\small $13$&\small $23$&\small $37$&\small $58$&\small $87$\\
&\small $1$&\small $3$&\small $7$&\small $13$&\small $22$&\small $34$&\small $50$&\small $70$\\
\hline
\multirow{4}{*}{$3$}&&\small $1$&\small $4$&\small $13$&\small $38$&\small $108$&\small $325$&\small $1275$\\
&&\small $1$&\small $4$&\small $13$&\small $37$&\small $101$&\small $284$&\small $956$\\
&&\small $1$&\small $4$&\small $13$&\small $37$&\small $101$&\small $284$&\small $956$\\
&&\small $1$&\small $4$&\small $13$&\small $37$&\small $92$&\small $209$&\small $442$\\
\hline
\multirow{4}{*}{$4$}&&&\small $1$&\small $5$&\small $23$&\small $108$&\small $940$&\small $190214$\\
&&&\small $1$&\small $5$&\small $23$&\small $101$&\small $677$&\small $70569$\\&&&\small $1$&\small $5$&\small $23$&\small $101$&\small $644$&\small $55081$\\&&&\small $1$&\small $5$&\small $23$&\small $100$&\small $432$&\small $1804$\\
\hline
\multirow{4}{*}{$5$}&&&&\small $1$&\small $6$&\small $37$&\small $325$&\small $190214$\\
&&&&\small $1$&\small $6$&\small $37$&\small $284$&\small $70569$\\
&&&&\small $1$&\small $6$&\small $37$&\small $284$&\small $55081$\\
&&&&\small $1$&\small $6$&\small $37$&\small $272$&\small $2806$\\
\hline
\multirow{4}{*}{$6$}&&&&&\small $1$&\small $7$&\small $58$&\small $1275$\\
&&&&&\small $1$&\small $7$&\small $58$&\small $956$\\
&&&&&\small $1$&\small $7$&\small $58$&\small $956$\\
&&&&&\small $1$&\small $7$&\small $58$&\small $817$\\
\hline

\end{tabular}
\caption{Numbers of matroids, base-orderable matroids, strongly base-orderable matroids and transversal matroids.}
\label{tab:transversal}
\end{center}
\end{table}

\section{Matroids on ten elements?}

Given that 30+ years have elapsed since the catalogue of matroids on 8 elements was created and with the benefit of advances both in raw computational power and techniques in combinatorial construction, it may seem rather unambitious to extend the catalogue only to 9 elements. 

However our initial experiments on the feasibility of constructing the matroids on 10 elements lead us to the conclusion that even {\em counting} the 10-element matroids would be a major undertaking, let alone constructing them. 

This is very unfortunate because we have a strong feeling that the absence of geometric representations available for rank-3 and rank-4 matroids means that rank-5 matroids are in some sense much less well understood than their lower rank counterparts. One of our original motivations
in embarking on this project was the belief that the rank-5 matroids on 10 elements might be a fertile source of interesting and/or counterintuitive examples and counterexamples for this reason.

\subsection{Paving matroids of rank 4}

From our analysis above, the sparse paving matroids of rank 4 on 10 elements are in 1-1 correspondence with independent sets in the Johnson graph $J(10,4)$, with isomorphism of
matroids and isomorphism under the automorphism group $S_{10}$ of the graph being
the same. Therefore a straightforward orderly algorithm as outlined in Royle \cite{MR1614301} can be used to construct them. This computation was performed in a few days using idle time on a network of about 50 computers, and the resulting numbers are presented 
in Table~\ref{tab:spm10} which shows a total of
$3150333219$ (i.e. $\approx 3.150 \times 10^9$) sparse paving matroids of rank 4 on 10 elements.

\begin{table}
\begin{center}
\begin{tabular}{|cr|cr|cr|cr|}
\hline
Size&Number&Size&Number&Size&Number&Size&Number\\
\hline
$0$&$1$&$8$&$521367$&$16$&$579539500$&$24$&$1355$\\
$1$&$2$&$9$&$3539486$&$17$&$329728133$&$25$&$250$\\
$2$&$3$&$10$&$18146294$&$18$&$130254690$&$26$&$58$\\
$3$&$13$&$11$&$69516384$&$19$&$35087875$&$27$&$13$\\
$4$&$73$&$12$&$197898106$&$20$&$6400127$&$28$&$4$\\
$5$&$575$&$13$&$416277780$&$21$&$818999$&$29$&$1$\\
$6$&$5838$&$14$&$642315652$&$22$&$84722$&$30$&$1$\\
$7$&$59818$&$15$&$720126836$&$23$&$9263$&&\\
\hline
\end{tabular}
\end{center}
\caption{Independent sets in $J(10,4)$}
\label{tab:spm10}
\end{table}

Computation of the {\em non-sparse} paving matroids of rank 4 on 10 elements is a somewhat fiddly bookkeeping exercise, but it involves
no qualitatively different techniques. The essence of our approach is to 
divide the search according to 
whether the largest hyperplane
has size $k=5$, $6$, $7$, $8$ or $9$.
For each size $k$, we construct an auxiliary graph $G(k)$ defined on the $4$-, $5$-, $\ldots$, $k$-sets that meet a fixed $k$-set (e.g. $\{0,1,\ldots,k-1\}$) in less than 3 points and with adjacency again defined by intersection in at least 3 points. Then an independent set of $G(k)$ together with $\{0,1,\ldots,k-1\}$ and all 3-sets not already covered forms the set of hyperplanes of a non-sparse paving matroid. However we need to be a little careful with isomorphism --- this procedure distinguishes a particular $k$-set and so if a matroid has $c$ orbits on hyperplanes of size $k$, then it will contribute $c$
pairwise non-isomorphic independent sets to $G(k)$. Therefore each independent set 
contributes $1/c$ to the
total count of matroids, where $c$ is the number of orbits that the corresponding matroid has
on hyperplanes of size $k$. Of course if the matroid has only one hyperplane of size $k$,
then $c=1$ follows immediately with no special calculation.

Table~\ref{tab:nspm10} shows the results of this calculation broken down according to the size of the
largest hyperplane $k$ and how many hyperplanes of this size are in the matroid.

\begin{table}
\begin{center}
\begin{tabular}{|ccr|}
\hline
Max. hyp. $k$ & No. $k$-hyps. & No. matroids\\
\hline
5&1&1222076172\\
5&2&147724716\\
5&3&5558695\\
5&4&64194\\
5&5& 232\\
5&6&6\\
6&1& 2369590\\
6&2& 164\\
7&1& 435\\
8&1&5\\
9&1&1\\
\hline
Total & &1377794210\\
\hline
\end{tabular}
\end{center}
\caption{Non-sparse paving matroids of rank 4 on 10 elements}
\label{tab:nspm10}
\end{table}

Adding the numbers of sparse and non-sparse paving matroids, we conclude that there
are  $4528127429$ ($\approx 4.528 \times 10^9$) paving matroids of rank 4 on 10 elements.

\subsection{Paving matroids of rank 5}

We have been unable to complete the analogous computation for the sparse paving matroids of
rank 5 on 10 elements. These correspond to independent sets in the Johnson graph $J(10,5)$ but
with one additional complication. The automorphism group of $J(10,5)$ is $S_{10} \times Z_2$ 
with the additional $Z_2$ being induced by complementation of 5-sets. This
means that each independent set of $J(10,5)$ produced by the orderly algorithm corresponds to 
a dual pair of matroids --- usually two matroids, but only one when the matroid is self-dual. 
Thus the total number of matroids is twice the number of independent sets of $J(10,5)$ minus the number of self-dual matroids.

We can determine the number of self-dual sparse paving matroids on 10 elements in a separate computation by exploiting the fact that the corresponding independent sets must have a non-trivial automorphism involving the $Z_2$ part of the automorphism group of $J(10,5)$. This separate computation yields a total of $99022169$ self-dual sparse paving matroids. 

However the sheer number of independent sets in $J(10,5)$ makes it infeasible for us to complete the first part of the computation. We can however make an ``informed guess'' of the magnitude of the number by executing a fixed percentage of the search. First, the orderly algorithm was used to compute the entire collection of independent sets of size 9, of which there are 20680075. A random sample of this collection was selected, and the search completed just using these as the starting points. 
Although it is hard to say anything statistically precise, our prior experience with such orderly algorithms suggest that the number of independent sets produced is roughly proportional to the size of the random sample of starting points. 

A sample of 60000 starting points (0.2901\% of the search space) yielded $3.875 \times 10^9$ independent sets giving an estimate of $1.336 \times 10^{12}$ independent sets in $J(10,5)$.  Therefore we estimate that there are about $2.65 \times 10^{12}$ sparse paving matroids of rank 5 on 10 elements.

\bibliographystyle{acm2url}
\bibliography{masterbib}

\begin{thebibliography}{10}

\bibitem{MR740223}
{ Acketa, D.}
\newblock { The catalogue of all nonisomorphic matroids on at most 8 elements},
  vol.~1 of { Special Issue}.
\newblock University of Novi Sad Institute of Mathematics Faculty of Science,
  Novi Sad, 1983.

\bibitem{MR0419270}
{ Blackburn, J.~E., Crapo, H.~H., and Higgs, D.~A.}
\newblock A catalogue of combinatorial geometries.
\newblock { Math. Comp. 27\/} (1973), 155--166; addendum, ibid. 27 (1973), no.
  121, loose microfiche suppl. A12--G12.

\bibitem{MR0304210}
{ Brualdi, R.~A., and Dinolt, G.~W.}
\newblock Characterizations of transversal matroids and their presentations.
\newblock { J. Combinatorial Theory Ser. B 12\/} (1972), 268--286.

\bibitem{MR0190045}
{ Crapo, H.~H.}
\newblock Single-element extensions of matroids.
\newblock { J. Res. Nat. Bur. Standards Sect. B 69B\/} (1965), 55--65.

\bibitem{MR2119297}
{ Dukes, W. M.~B.}
\newblock On the number of matroids on a finite set.
\newblock { S\'em. Lothar. Combin. 51\/} (2004/05), Art. B51g, 12 pp.
  (electronic).

\bibitem{MR635936}
{ McKay, B.~D.}
\newblock Practical graph isomorphism.
\newblock In { Proceedings of the Tenth Manitoba Conference on Numerical
  Mathematics and Computing, Vol. I (Winnipeg, Man., 1980)\/} (1981), vol.~30,
  pp.~45--87.
\newblock Available from: \url{http://cs.anu.edu.au/~bdm/nauty/PGI/}.

\bibitem{MR1606516}
{ McKay, B.~D.}
\newblock Isomorph-free exhaustive generation.
\newblock { J. Algorithms 26}, 2 (1998), 306--324.

\bibitem{MR2065730}
{ Oxley, J.}
\newblock What is a matroid?
\newblock { Cubo Mat. Educ. 5}, 3 (2003), 179--218.

\bibitem{MR1207587}
{ Oxley, J.~G.}
\newblock { Matroid theory}.
\newblock Oxford Science Publications. The Clarendon Press Oxford University
  Press, New York, 1992.

\bibitem{MR0491273}
{ Read, R.~C.}
\newblock Every one a winner or how to avoid isomorphism search when
  cataloguing combinatorial configurations.
\newblock { Ann. Discrete Math. 2\/} (1978), 107--120.
\newblock Algorithmic aspects of combinatorics (Conf., Vancouver Island, B.C.,
  1976).

\bibitem{MR1614301}
{ Royle, G.~F.}
\newblock An orderly algorithm and some applications in finite geometry.
\newblock { Discrete Math. 185}, 1-3 (1998), 105--115.

\bibitem{oeis}
{ Sloane, N. J.~A.}
\newblock The {O}n-{L}ine {E}ncylopaedia of {I}nteger {S}equences.
\newblock Available from: \url{http://www.research.att.com/~njas/sequences}.

\bibitem{MR0278975}
{ Welsh, D. J.~A.}
\newblock Combinatorial problems in matroid theory.
\newblock In { Combinatorial Mathematics and its Applications (Proc. Conf.,
  Oxford, 1969)}. Academic Press, London, 1971, pp.~291--306.

\bibitem{MR0427112}
{ Welsh, D. J.~A.}
\newblock { Matroid theory}.
\newblock Academic Press [Harcourt Brace Jovanovich Publishers], London, 1976.
\newblock L. M. S. Monographs, No. 8.

\bibitem{MR2179649}
{ Whittle, G.}
\newblock Recent work in matroid representation theory.
\newblock { Discrete Math. 302}, 1-3 (2005), 285--296.

\end{thebibliography}

 \end{document}